\documentstyle[12pt]{article}
\newcommand{\const}{\mathop{\rm const}\limits}

\newcommand{\supp}{\mathop{\rm supp}\limits}

\newcommand{\diam}{\mathop{\rm diam}\limits}

\newcommand{\mes}{\mathop{\rm mes}\limits}

\newcommand{\card}{\mathop{\rm card}\limits}

\begin{document}

\begin{center}

{\bf  CONTINUITY OF FUNCTIONS BELONGING TO THE} \par

\vspace{4mm}

{\bf FRACTIONAL ORDER  SOBOLEV-GRAND LEBESGUE SPACES}\par

\vspace{4mm}

 $ {\bf E.Ostrovsky^a, \ \ L.Sirota^b } $ \\

\vspace{4mm}

$ ^a $ Corresponding Author. Department of Mathematics and computer science, Bar-Ilan University, 84105, Ramat Gan, Israel.\\
\end{center}
E - mail: \ galo@list.ru \  eugostrovsky@list.ru\\
\begin{center}
$ ^b $  Department of Mathematics and computer science. Bar-Ilan University,
84105, Ramat Gan, Israel.\\

E - mail: \ sirota3@bezeqint.net\\

\vspace{3mm}
                    {\sc Abstract.}\\

 \end{center}

 \vspace{4mm}

 We extend in this  article the classical  Sobolev’s inequalities for the
module of continuity for the functions belonging to the integer order Sobolev's space on
the {\it fractional order} Sobolev - Bilateral Grand Lebesgue spaces.\par
 As a consequence, we deduce the fractional Orlicz - Sobolev imbedding theorems
and investigate the rectangle module of continuity  of non-Gaussian  multiparameter random fields. \par
  \vspace{4mm}


\vspace{4mm}

 \par

\vspace{4mm}

\section{Notations. Statement of problem.}

\vspace{3mm}

 Let $ D $ be convex non-empty bounded closed domain with Lipschitz boundary in the whole space $ R^d, \ d = 1,2,\ldots, $
and let $ f: D \to R $ be measurable function. We assume further for simplicity that $ D = [0,1]^d. $
 We denote $ |x| = (\sum_{i=1}^d x_i^2 )^{1/2},  \ \alpha  = \const \in (0,1], $

$$
|f|_p = |f|_{p,D} = \left[ \int_D |f(x)|^p \ dx  \right]^{1/p},  \ |u(\cdot, \cdot)|_p = |u(\cdot, \cdot)|_{p,D^2} =
$$

$$
 \left[ \int_D \int_D |u(x,y)|^p \ dx dy \right]^{1/p}, \ p = \const \ge 1,
$$

$$
  \omega(f,\delta) = \sup \{|f(x) - f(y)|: \ x,y \in D,  |x-y| \le \delta  \},  \ \delta  \in [0, \diam(D)], \eqno(1.0)
$$

$$
G_{\alpha}[f] (x,y)= \frac{f(x) - f(y) }{|x-y|^{\alpha}},  \hspace{5mm}  \nu(dx,dy) = \frac{dx dy}{|x-y|}, \eqno(1.1)
$$

$$
|u(\cdot, \cdot)|_{p,\nu} = |u(\cdot, \cdot)|_{p, \nu, D^2} = \left[ \int_D \int_D |u(x,y)|^p \ \nu(dx, dy) \right]^{1/p}, \eqno(1.2)
$$

$$
||f||W(\alpha,p) = | G_{\alpha}[f] (\cdot, \cdot) |_{p, \nu, D^2}. \eqno(1.3)
$$
 The norm  $ ||\cdot||W(\alpha,p), $ more precisely, semi-norm  is said to be {\it fractional } Sobolev's norm or similar
{\it Aronszajn, Gagliardo or Slobodeckij } norm; see, e.g. \cite{Nezzaa1}.\par

  If in the definition (1.3) instead the $ L_p(D^2)  $ stands another norm $ ||\cdot|| V(D^2), $ for instance,  Lorentz,
 Marcinkiewicz or Grand Lebesgue, (we recall its definition further),   we obtain correspondingly the definition of the
 fractional  $ ||\cdot|| V(D^2) $ norm. \par

 The inequality

 $$
 |f(t) - f(s)| \le  8 \cdot 4^{1/p} \cdot \left[\frac{\alpha + 1/p}{\alpha - 1/p}\right]   \cdot |t-s|^{\alpha - 1/p} \cdot
 ||f||W(\alpha,p), \eqno(1.4)
 $$
or equally

$$
\omega(f,\delta) \le 8 \cdot 4^{1/p} \cdot \left[\frac{\alpha + 1/p}{\alpha - 1/p}\right]   \cdot \delta^{\alpha - 1/p} \cdot
\left[ \int_D \int_D \frac{|f(x) - f(y)|^p \ dx dy }{|x-y|^{\alpha p + 1}}   \right]^{1/p}, \eqno(1.5)
$$
which is true in the case $  d = 1 $ (the multidimensional case will be consider further), $ p > 1/\alpha, $ is called {\it fractional}
Sobolev, or Aronszajn, Gagliardo, Slobodeckij inequality.\par
 More precisely, the inequality (1.4) implies that the function $ f $ may be redefined  on the set of measure zero as a
continuous function  for which (1.4) there holds.\par
 Another look on the inequality (1.4): it may be construed as an imbedding theorem from the Sobolev fractional space into the
space of (uniform) continuous functions on the set $  D. $ \par
 The proof of the our version of inequality (1.4) may be obtained immediately from an article \cite{Hu1}, which based in turn on the
famous Garsia - Rodemich - Rumsey inequality, see \cite{Garsia1}.\par
 There are many generalizations of fractional Sobolev's imbedding theorem: on the Sobolev - Orlicz's spaces \cite{Adams1}, p. 253-364, on the
so-called {\it integer} Sobolev - Grand Lebesgue spaces  \cite{Ostrovsky100}, on the Lorentz and Marcinkiewicz spaces etc.\par

\vspace{3mm}
 {\bf  Our goal is to extend the  Sobolev's imbedding theorem from  integer Sobolev Grand Lebesgue spaces on the fractional
 Sobolev Grand Lebesgue spaces.}\par
\vspace{3mm}
 We recall further the definition of these spaces.\par
The applications of fractional Sobolev  and Sobolev-Grand Lebesgue spaces in the theory of Partial Differential Equations
are described, e.g. in \cite{Lieb1},   \cite{Nezzaa1}, \cite{Runst1}; in the Functional Analysis - in  \cite{Adams1}, \cite{Garsia1},
\cite{Kaminska1}, \cite{Lieb1}, \cite{Milman3};  in the theory of Random Processes and Fields - in \cite{Garsia1}, \cite{Hu1},
\cite{Inahama1}, \cite{Ral'chenko1}; see also reference therein. \par

\vspace{3mm}

We recall in the  remainder part of this section briefly  the definition of the so-called Grand Lebesgue spaces;   more detail
investigation of these spaces see in \cite{Fiorenza3}, \cite{Iwaniec2}, \cite{Kozachenko1}, \cite{Liflyand1}, \cite{Ostrovsky1},
\cite{Ostrovsky2}; see also reference therein.\par

  Recently  appear the so-called Grand Lebesgue Spaces $ GLS = G(\psi) =G\psi =
 G(\psi; A,B), \ A,B = \const, A \ge 1, A < B \le \infty, $ spaces consisting
 on all the measurable functions $ f: X \to R $ with finite norms

$$
   ||f||G(\psi) \stackrel{def}{=} \sup_{p \in (A,B)} \left[ |f|_p /\psi(p) \right]. \eqno(1.6)
$$

  Here $ \psi(\cdot) $ is some continuous positive on the {\it open} interval
$ (A,B) $ function such that

$$
     \inf_{p \in (A,B)} \psi(p) > 0, \ \psi(p) = \infty, \ p \notin (A,B).
$$
 We will denote
$$
 \supp (\psi) \stackrel{def}{=} (A,B) = \{p: \psi(p) < \infty, \}
$$

The set of all $ \psi $  functions with support $ \supp (\psi)= (A,B) $ will be
denoted by $ \Psi(A,B). $ \par
  This spaces are rearrangement invariant, see \cite{Bennet1}, and
    are used, for example, in the theory of probability  \cite{Kozachenko1},
  \cite{Ostrovsky1}, \cite{Ostrovsky2}; theory of Partial Differential Equations \cite{Fiorenza3},
  \cite{Iwaniec2};  functional analysis \cite{Fiorenza3}, \cite{Iwaniec2},  \cite{Liflyand1},
  \cite{Ostrovsky2}; theory of Fourier series, theory of martingales, mathematical statistics,
  theory of approximation etc.\par

 Notice that in the case when $ \psi(\cdot) \in \Psi(A,\infty)  $ and a function
 $ p \to p \cdot \log \psi(p) $ is convex,  then the space
$ G\psi $ coincides with some {\it exponential} Orlicz space. \par
 Conversely, if $ B < \infty, $ then the space $ G\psi(A,B) $ does  not coincides with
 the classical rearrangement invariant spaces: Orlicz, Lorentz, Marcinkiewicz  etc.\par

  The fundamental function of these spaces $ \phi(G(\psi), \delta) = ||I_A ||G(\psi), \mes(A) = \delta, \ \delta > 0, $
where  $ I_A  $ denotes as ordinary the indicator function of the measurable set $ A, $ by the formulae

$$
\phi(G(\psi), \delta) = \sup_{ p \in \supp (\psi)} \left[ \frac{\delta^{1/p}}{\psi(p)} \right].  \eqno(1.7)
$$
The fundamental function of arbitrary rearrangement invariant spaces plays very important role in functional analysis,
theory of Fourier series and transform \cite{Bennet1} as well as in our further narration. \par

 Many examples of fundamental functions for some $ G\psi $ spaces are calculated in  \cite{Ostrovsky1}, \cite{Ostrovsky2}.\par

\vspace{3mm}

{\bf Remark 1.1} If we introduce the {\it discontinuous} function

$$
\psi_{(r)}(p) = 1, \ p = r; \psi_{(r)}(p) = \infty, \ p \ne r, \ p,r \in (A,B)
$$
and define formally  $ C/\infty = 0, \ C = \const \in R^1, $ then  the norm
in the space $ G(\psi_r) $ coincides with the $ L_r $ norm:

$$
||f||G(\psi_{(r)}) = |f|_r.
$$
Thus, the Grand Lebesgue Spaces are direct generalization of the
classical exponential Orlicz's spaces and Lebesgue spaces $ L_r. $ \par

\vspace{3mm}

{\bf Remark 1.2}  The function $ \psi(\cdot) $ may be generated as follows. Let $ \xi = \xi(x)$
be some measurable function: $ \xi: X \to R $ such that $ \exists  (A,B):
1 \le A < B \le \infty, \ \forall p \in (A,B) \ |\xi|_p < \infty. $ Then we can
choose

$$
\psi(p) = \psi_{\xi}(p) = |\xi|_p.
$$

 Analogously let $ \xi(t,\cdot) = \xi(t,x), t \in T, \ T $ is arbitrary set,
be some {\it family } $ F = \{ \xi(t, \cdot) \} $ of the measurable functions:
$ \forall t \in T  \ \xi(t,\cdot): X \to R $ such that
$$
 \exists  (A,B): 1 \le A < B \le \infty, \ \sup_{t \in T} \
|\xi(t, \cdot)|_p < \infty.
$$
 Then we can choose

$$
\psi(p) = \psi_{F}(p) = \sup_{t \in T}|\xi(t,\cdot)|_p.
$$
The function $ \psi_F(p) $ may be called as a {\it natural function} for the family $ F. $
This method was used in the probability theory, more exactly, in
the theory of random fields, see \cite{Kozachenko1},\cite{Ostrovsky1}, chapters 3,4. \par

\vspace{3mm}

{\bf Remark 1.3} Note that the so-called {\it exponential} Orlicz spaces are particular cases of
Grand Lebesgue spaces  \cite{Kozachenko1}, \cite{Ostrovsky1}, p. 34-37.  In detail, let the $ N- $
Young-Orlicz function has a view

$$
 N(u) = e^{\mu(u)}, \eqno(1.8)
$$
where the function $ u \to \mu(u) $ is convex even twice differentiable function such that

$$
\lim_{u \to \infty} \mu'(u) = \infty.
$$
 Introduce a new function
$$
\psi_{\{N\}}(x) = \exp \left\{ \frac{\left[\log N(e^x) \right]^*}{x}   \right\}, \eqno(1.9)
$$
where $  g^*(\cdot) $ denotes the Young-Fenchel transform of the function $  g:  $

$$
g^*(x) = \sup_y (xy - g(y)).
$$
 Conversely,  the $  N  - $ function may be calculated up to equivalence
  through corresponding function $ \psi(\cdot) $  as follows:

 $$
 N(u) = e^{\tilde{\psi}^*(\log |u|) }, \ |u| > 3; \ N(u) = C u^2, |u| \le 3; \  \tilde{\psi}(p) = p \log \psi(p). \eqno(1.10)
 $$
 The Orlicz's space $ L(N) $ over our probabilistic space is equivalent up to sublinear norms equality with
Grand Lebesgue space $ G\psi_{\{N\}}. $ \par

 For instance, if $ N(u) = N_2(u) := \exp(u^2/2) - 1, $ then $ \psi_{\{N_2\}}(p) \asymp \sqrt{p}, \ p \ge 1. $
The centered r.v.  belonging  to the Orlicz's space $ L(N_2)$ are called subgaussian. \par
  More generally, if $ N(u) = N_m(u) := \exp(u^m/m) - 1, \ m = \const > 0, $ then $ \psi_{\{N_m\}}(p) \asymp p^{1/m}, \ p \ge 1. $ \par

\vspace{3mm}

{\bf Remark 1.4.} The theory of probabilistic {\it exponential} Grand Lebesgue spaces
or equally exponential Orlicz spaces gives a
  very convenient apparatus for investigation of
the r.v. with exponential decreasing tails of distributions. Namely, the non-zero  r.v. $ \eta $ belongs to the
Orlicz space $ L(N), $  where $ N = N(u) $ is function described in equality (1.8), if and only if

$$
{\bf P} (\max(\eta, -\eta) > z) \le \exp(-\mu(C z)), \ z > 1,  \ C = C(N(\cdot), ||\eta||L(N)) \in (0,\infty). \eqno(1.11)
$$
(Orlicz's version). \par
 Analogously may be written a Grand Lebesgue version of this inequality.
 In detail,  if $ 0 < ||\eta||G\psi< \infty,  $ then

 $$
{\bf P} (\max(\eta, -\eta) > z) \le  2\exp \left(- \tilde{\psi}(\log [z /||\eta||G\psi] ) \right), z \ge ||\eta||G\psi.
 \eqno(1.12)
$$
  Conversely, if

 $$
{\bf P} (\max(\eta, -\eta) > z) \le  2\exp \left(- \tilde{\psi}(\log [z /K] ) \right), z \ge K,
$$
then $ ||\eta||G\psi \le C(\psi)  \cdot K, \ C(\psi) \in (0,\infty). $  \par

\vspace{3mm}
\section{One dimensional result.}

\vspace{3mm}

  Let as before $ \alpha = \const \in (0,1]; $  the case $  \alpha > 1  $ in this section is trivial.   We introduce the  $ \Psi $ function
$ \zeta_{\alpha}(p) $ as follows:

$$
\zeta_{\alpha}(p):=||f||W(\alpha,p), \  (A,B):= \supp \ [\zeta_{\alpha}(\cdot)] \eqno(2.0)
$$
 and suppose $  1 \le A < B \le \infty.  $ \par

  Denote $ A(\alpha) = \max(A, 1/\alpha) $ and suppose also
$$
A(\alpha) < B. \eqno(2.1)
$$
 Obviously, the restriction (2.1) is satisfied always in the case $  B = \infty. $ \par
 We define a new  psi - function $ \psi_{\alpha}(p)  $ as follows.

$$
\psi_{\alpha}(p):= \zeta_{\alpha}(p) \cdot 8 \cdot 4^{1/p} \cdot \frac{ \alpha + 1/p}{\alpha - 1/p}, \ p \in (A(\alpha), B) \eqno(2.2)
$$
and $ \psi_{\alpha}(p) = \infty  $ otherwise. \par

 The fractional  Sobolev - Grand Lebesgue norm $ ||f||S(\alpha, \psi) $ for arbitrary function $ \psi \in \Psi $
  of the function $ f: D \to R $ may be defined
in accordance with first section up to multiplicative constant as follows: $ ||f||S(\alpha,\psi) \stackrel{def}{=}  $

$$
 \sup_{p \in (A(\alpha), B)} \left\{ \left[ 8 \cdot 4^{1/p} \cdot \left(\frac{\alpha + 1/p}{\alpha - 1/p}\right)  \cdot
 \left( \int_D \int_D \frac{|f(x) - f(y)|^p \ dx dy }{|x-y|^{\alpha p + 1}}   \right)^{1/p}\right] /\psi(p) \right\}, \eqno(2.3)
$$
so that  the function $ \psi_{\alpha}(p) $ is the natural function for the function $ G_{\alpha}(x,y) $  relative the two dimensional
measure $  \nu. $ \par

\vspace{3mm}

{\bf Theorem 2.1.} {\it  Let }  $  d = 1 $ {\it  and let the condition 2.1 be satisfied. Then}

$$
\omega(f, \delta) \le \frac{\delta^{\alpha}}{\phi(G\psi_{\alpha}, \delta)} \cdot  ||f||G \psi_{\alpha},
\hspace{5mm} \delta \in (0,\diam D). \eqno(2.4)
$$

\vspace{3mm}

{\bf Proof.} We can and will suppose without loss of generality $ D = [0,1]  $  and $ ||f||S(\alpha, \psi_{\alpha}) = 1. $
 Let $ p \in (A(\alpha), B).  $ It follows from the definition  for $ ||f||S(\alpha,\psi_{\alpha}) $ (2.3)  that

 $$
 8 \cdot 4^{1/p} \cdot \left[\frac{\alpha + 1/p}{\alpha - 1/p}\right]  \cdot
\left[ \int_D \int_D \frac{|f(x) - f(y)|^p \ dx dy }{|x-y|^{\alpha p + 1}}   \right]^{1/p} \le \psi_{\alpha}(p). \eqno(2.5)
$$
 The application of estimate (1.5) yields

 $$
 \omega(f,\delta) \le \delta^{\alpha - 1/p} \ \psi_{\alpha}(p),
 $$
or equally

$$
\frac{\omega(f,\delta)}{\delta^{\alpha}} \le \frac{1}{ \delta^{1/p}/\psi_{\alpha}(p) }.
$$
 Since the value $ p  $ is arbitrary in the interval $ p \in (A(\alpha), B),  $ we conclude

$$
\frac{\omega(f,\delta)}{\delta^{\alpha}} \le  \inf_{p \in (A(\alpha), B)} \left[ \frac{1}{ \delta^{1/p}/\psi_{\alpha}(p) } \right] =
$$
$$
\frac{1}{ \sup_{p \in (A(\alpha), B)} [\delta^{1/p}/\psi_{\alpha}(p)] } = \frac{1}{\phi(G(\psi_{\alpha}),\delta)} =
 \frac{ ||f||S(\alpha, \psi_{\alpha})}{\phi(G(\psi_{\alpha}),\delta)}, \eqno(2.6)
$$
Q.E.D. \par

\vspace{3mm}

{\bf Example 2.1.}  Suppose
$$
\psi_{\alpha}(p) \asymp \psi^{(a,b; A,B)}(p) := (p-A)^{-a} (B-p)^{-b}, \ p \in (A,B),  \ A \ge 1/\alpha, \ a,b = \const \in (0,\infty).
$$
 The fundamental function for the spaces  $ G\psi^{(a,b; A,B)} $ is investigated in \cite{Ostrovsky2}. Take note only that as $ \delta \to 0+ $

$$
\phi(G\psi^{(a,b; A,B)}, \delta ) \asymp \delta^{1/B} \ |\log \delta|^{-b}, \ 0 < \delta < 1/e.
$$
 Therefore in the considered case

 $$
 \omega(f,\delta) \le C(\alpha,a,b; A,B) \cdot \delta^{\alpha - 1/B} \cdot |\log \delta|^b \cdot ||f||(G\psi^{(a,b; A,B)},  \ 0 < \delta < 1/e.
 $$

\vspace{3mm}

{\bf Example 2.2.} Let now

$$
\psi_{\alpha}(p) \asymp \psi_{[\beta]}(p) := p^{\beta}, \ \beta = \const > 0, p > 1/\alpha.
$$
We find analogously the example 2.1
$$
\phi(G\psi_{[\beta]},\delta) \asymp |\log \delta|^{-\beta}, \
\omega(f,\delta) \le C(\alpha,\beta) \cdot  \delta^{\alpha} \ |\log \delta|^\beta  \cdot ||f||G\psi_{[\beta]},  \ 0 < \delta < 1/e.
$$

\vspace{3mm}

{\bf Remark 2.1.} Assume in addition to the condition of theorem 2.1 $ \alpha = 1. $ Recall that for arbitrary rearrangement
invariant space $ X $

$$
\phi(X',\delta) = \frac{\delta}{\phi(X,\delta)},
$$
see \cite{Bennet1}, chapter 3; here $ X' $ denotes the associate space to the space $ X. $ \par

 The conclusion of theorem 2.1 in the considered case $ \alpha = 1 $ may be rewritten as follows:

$$
\omega(f, \delta) \le  \phi((G\psi)',\delta) \cdot ||f||S(1, \psi_1), \hspace{5mm} \delta \in (0,\diam D).  \eqno(2.7)
$$
  In general case $ \alpha  \in (0,1] $ we have

$$
\omega(f, \delta) \le  \delta^{\alpha - 1} \cdot \phi((G\psi_{\alpha})',\delta) \cdot ||f||S(\alpha, \psi_{\alpha}), \hspace{5mm} \delta \in (0,\diam D).  \eqno(2.8)
$$

\vspace{3mm}

{\bf  Remark 2.2.} Instead the function $ \nu_{\alpha} (p) $ may be used arbitrary it majorant.\par

\vspace{3mm}
{\bf  Remark 2.3.}  We discus the exactness of the assertion of theorem 2.1  further; now  we proceed to the consideration of
 multidimensional case. \par

\vspace{3mm}

{\bf  Remark 2.4.} If we use in the capacity  of the function $ \nu_{\alpha} (p) $ the discontinuous function  $ \psi_{(r)}(p),  $
we return  to the inequality 1.5.\par

\vspace{3mm}

{\bf  Remark 2.5.} If the value $  \alpha $ in the assertion (2.4) of theorem 2.1 is {\it variable} in some interval
$ \alpha \in (\alpha_-, \alpha_+), $ then obviously

$$
\omega(f, \delta) \le
\inf_{ \alpha \in (\alpha_-, \alpha_+)  } \left[ \frac{\delta^{\alpha}}{\phi(G\psi_{\alpha}, \delta)} \cdot  ||f||G \psi_{\alpha} \right],
\hspace{5mm} \delta \in (0,\diam D). \eqno(2.9)
$$

\vspace{3mm}
\section{Multi-dimensional result.}

\vspace{3mm}

The multidimensional case $ d= \dim D = 2,3,\ldots $  is more complicated. Suppose for simplicity $ D = [0.1]^d. $
This imply that $ x \in D \Leftrightarrow x = \vec{x} = (x_1,x_2, \ldots, x_d), \ 0 \le x_i \le 1. $\par

 We define as in \cite{Ral'chenko1}, \cite{Hu1} the {\it rectangle difference}  operator $ \Box[f](\vec{x}, \vec{y}) =
 \Box[f](x,y),  \ x,y \in D, \ f:D \to R $ as follows.

$$
\Delta^{(i)}[f](x,y) := f(x_1,x_2, \ldots, x_{i-1}, y_i, x_{i+1}, \ldots, x_d) - f(x_1,x_2, \ldots, x_{i-1}, x_i, x_{i+1}, \ldots, x_d),
$$
with obvious modification when $ i=1 $ or $ i=d; $
$$
\Box[f](x,y) \stackrel{def}{=} \left\{ \otimes _{i=1}^d \Delta^{(i)} \right\} [f](x,y).\eqno(3.1)
$$
 For instance, if $ d=2, $ then

$$
\Box[f](x,y) = f(y_1,y_2) - f(x_1,y_2) - f(y_1,x_2) + f(x_1,x_2).
$$

 If the function $  f: [0,1]^d \to R $ is $ d $ times continuous differentiable, then

 $$
 \Box[f](\vec{x},\vec{y}) = \int_{x_1}^{y_1}
 \int_{x_2}^{y_2} \ldots   \int_{x_d}^{y_d} \frac{ \partial^d f}{\partial x_1  \partial x_2  \ldots   \partial x_d } \ dx_1 dx_2 \ldots   dx_d .
 $$

 The {\it rectangle module of continuity  } $  \Omega(f, \vec{\delta} ) =  \Omega(f, \delta )  $  for the
(continuous  a.e.) function $ f $ and vector $ \vec{\delta} = \delta = ( \delta_1, \delta_2, \ldots, \delta_d)
\in [0,1]^d  $  may be defined as well as ordinary module of continuity $ \omega(f,\delta) $ as follows:

$$
\Omega(f, \vec{\delta} ) \stackrel{def}{=} \sup \{ |\Box[f](x,y)|, \ (x,y): |x_i - y_i| \le \delta_i, \ i = 1,2,\ldots,d \}.
$$

   Let $ \vec{\alpha} = \{  \alpha_k \}, \ \alpha_k \in (0,1], \ k=1,2,\ldots,d; \ p > p_0 \stackrel{def}{=}\max_k (1/\alpha_k),
  \ M = \card  \{i, \alpha_i = \min_k \alpha_k\}, \delta_i = |x_i - y_i|, \ \vec{\delta} = \{ \delta_i \}, i = 1,2,\ldots,d; $

$$
\vec{x}^{\vec{\alpha}} := \prod_{i=1}^d x_i^{\alpha_i}, \  \vec{\delta}^{\pm 1/p} :=  \left[\prod_{i=1}^d  \delta_i  \right]^{\pm 1/p},
$$

$$
G_{\vec{\alpha}}[f] (x,y) =  \frac{\Box[f](x,y) }{  |(\vec{x}- \vec{y})^{\vec{\alpha}} | },  \hspace{5mm}  \nu(dx,dy) = \frac{\vec{dx} \vec{ dy}}{|x-y|},
$$

$$
||f||W( \vec{\alpha},p) = | G_{\vec{\alpha}}[f] (\cdot, \cdot) |_{p, \nu, D^2}.
$$

 The norm  $ ||\cdot||W( \vec{\alpha},p), $ more precisely, semi-norm  is said to be {\it multidimensional fractional } Sobolev's norm or similar
{\it Aronszajn, Gagliardo or Slobodeckij } norm.  \par
 Define also as well as in the second section

$$
\zeta_{\vec{\alpha}}(p):=||f||W(\vec{\alpha},p), \  (A,B):= \supp \ [\zeta_{\vec{\alpha}}(\cdot)]
$$
 and suppose $  1 \le A < B \le \infty.  $ \par

  Denote $ A (\vec{\alpha}) = \max(A, p_0) $ and suppose also $ A(\vec{\alpha}) < B. $  \par

 We define a new  psi - function $ \psi_{\alpha}(p)  $ as follows.

$$
\psi_{\vec{\alpha}}(p):= \zeta_{\vec{\alpha}}(p) \cdot 8^d \cdot 4^{d/p} \cdot \prod_{k=1}^d \left[ \frac{ \alpha_k + 1/p}{\alpha_k - 1/p} \right], \
p \in (A(\vec{\alpha}), B).
$$

\vspace{3mm}

{\bf Theorem 3.1.} {\it  Let }  $  d = 1 $ {\it  and let all our condition be satisfied. Then
then there exists  a continuous modification on the set of measure zero
of the function $ f, $ which we will denote again $ f, $ for which}

$$
\Omega(f, \vec{ \delta}) \le \frac{\vec{\delta}^{\vec{\alpha}}}{\phi(G\psi_{\vec{\alpha} }, \prod_{k=1}^d \delta_k)} \cdot  ||f||G\psi_{\vec{\alpha}}.
$$

\vspace{3mm}

 {\bf Proof} is similar to one  in theorem 2.1.   \hspace{5mm} We can take as before $ ||f||G\psi_{\vec{\alpha}} = 1.$ \par

 The multidimensional  Garsia - Rodemich - Rumsey inequality was done  by
Konstantin Ral'chenko  \cite{Ral'chenko1} (2007) at $ d=2 $  and Yaozhong Hu and Khoa Le   \cite{Hu1} (2012)  in general case.
 Namely,  let $ \Phi = \Phi(y) $ be Young-Orlicz continuous even strictly increasing on the right-hand  semi-axis  function such that

 $$
 \Phi(0) = 0, \ \lim_{y \to \infty} \Phi(y) = \infty.
 $$
Let also $  p_k = p_k(u), \ u \in [0,1], \ k=1,2,\ldots,d $ be continuous strictly increasing  functions. We denote

$$
B := \int \int \ldots \int_{[0,1]^{2d}} \Phi \left[ \frac{|\Box [f](x,y)|}{\prod_{k=1}^d p_k(|x_k-y_k|)} \right] dx dy
$$
and assume $ B < \infty. $ Then

$$
|\Omega(f, \vec{ \delta}) | \le 8^d \int_0^{\delta_1} \int_0^{\delta_2} \ldots \int_0^{\delta_d}
\Phi^{-1}\left[ \frac{4^d B}{\prod_{j=1}^d u^2_j} \right] dp_1(u_1) dp_2(u_2) \ldots d p_d(u_d).
$$

 We refer further using for us the particular case of this inequality.
It asserts that for some modification of the multidimensional version of the function $  f $

$$
|\Box[f](x,y)| \le  8^d \ 4^{d/p} \ \prod_{i=1}^d \left[\frac{\alpha_i+1/p}{\alpha_i-1/p} \right] \cdot
\vec{\delta}^{\vec{\alpha}} \cdot  \vec{\delta}^{- 1/p} \times
$$

$$
 \left[ \int \int \ldots \int_{[0,1]^d} \frac{|\Box[f](x,y)|^p}{\prod_{k=1}^d |x_k - y_k|^{\alpha_k p + 1}}  dx dy \right]^{1/p} < \infty; \eqno(3.2)
$$
or equally

 $$
 \Omega(f,\vec{\delta}) \le \vec{\delta}^{\vec{\alpha}} \cdot \vec{\delta}^{ - 1/p} \cdot \psi_{\vec{\alpha}}(p),\eqno(3.3)
 $$
whence

$$
\frac{\Omega(f,\vec{\delta})}{\vec{\delta}^{\vec{\alpha}} }  \le \frac{1}{ \vec{\delta}^{1/p}/\psi_{\vec{\alpha}}(p) } =
\frac{1}{ [\prod_{k=1}^d \delta_k^{1/p}]/\psi_{\vec{\alpha}}(p) }. \eqno(3.4)
$$

 Since the value $ p  $ is arbitrary in the interval $ p \in (A(\vec{\alpha}), B),  $ we conclude

$$
\frac{\Omega(f, \vec{\delta})}{\delta^{\vec{\alpha}}} \le  \inf_{p \in (A(\vec{\alpha}), B)}
 \frac{1}{ [\prod_{k=1}^d \delta_k^{1/p} ]/ \psi_{\vec{\alpha}}(p) } =
$$
$$
\frac{1}{ \sup_{p \in (A(\vec{\alpha}), B)} [\prod_{k=1}^d \delta_k^{1/p}/\psi_{\vec{\alpha}}(p)] } =
 \frac{1}{\phi(G(\psi_{\vec{\alpha}}),\prod_{k=1}^d \delta_k)} =
 \frac{ ||f||G\psi_{\vec{\alpha}}}{\phi(G(\psi_{\vec{\alpha}}),\prod_{k=1}^d \delta_k)},
$$
Q.E.D. \par

\vspace{3mm}

\section{Application to the theory of random fields. }

\vspace{3mm}

  Let $ \xi = \xi(x) = \xi(x_1,x_2,\ldots, x_d) = \xi(\vec{x}), \ x_i \in [0,1] $ be separable random field (r.f),
  not necessary to be Gaussian.
   The correspondent probability and expectation we will denote  by $ {\bf P}, \ {\bf E,}  $  and  the probabilistic
 Lebesgue-Riesz $ L_p $   norm of a random variable (r.v) $  \eta $  we will denote as follows:

 $$
 |\eta|_p \stackrel{def}{=} \left[ {\bf E} |\eta|^p \right]^{1/p}.
 $$

   We find in this section some sufficient condition for continuity of $ \xi(x) $  and estimates for it {\it rectangle} modulus of continuity
 $ \Omega(f,\vec{\delta}). $  We apply in this section the results obtained before.
   Recall that the first publication about fractional  Sobolev's  inequalities
   \cite{Garsia1}  was devoted in particular to the such a problem; see also articles \cite{Hu1}, \cite{Ral'chenko1}. \par
\vspace{3mm}

 Let us introduce the following  natural $ \Psi $ function: $  \theta_{\vec{\alpha}}(p) = $

$$
\theta_{\alpha}(p) =  8^d \cdot 4^{d/p} \cdot \prod_{k=1}^d \left[ \frac{\alpha_k + 1/p}{\alpha_k - 1/p} \right] \cdot
\left[ \int_0^1 \int_0^1 {\bf E} |G_{\vec{\alpha}}[\xi](x,y)|^p \nu(dx,dy) \right]^{1/p}, \eqno(4.1)
$$

$$
\alpha = \vec{\alpha} = \{ \alpha_1, \alpha_2, \ldots,\alpha_d \}, \ \alpha_k = \const > 0;
$$
and suppose the function $ \theta_{\alpha}(p) $ has non-trivial support such that

$$
A = \inf \supp \theta_{\alpha}(\cdot) \ge 1/\min_k\alpha_k, \ B = \sup \supp \theta_{\alpha} \in (A, \infty].
$$

\vspace{3mm}

{\bf Theorem 4.1.}

$$
|\Omega[\xi],\vec{\delta}|_A  \le \frac{\vec{\delta}^{\vec{\alpha}}}{\phi(G\theta_{\alpha}, \prod_{k=1}^d \delta_k) }. \eqno(4.2)
$$
{\bf Proof.} We use the inequality (3.2):

$$
\Box[\xi](x,y)| \le  8^d \ 4^{d/p} \ \prod_{i=1}^d \left[\frac{\alpha_i+1/p}{\alpha_i-1/p} \right] \cdot
\vec{\delta}^{\vec{\alpha}} \cdot  \vec{\delta}^{- 1/p} \times
$$

$$
 \left[ \int \int \ldots \int_{[0,1]^d} \frac{|\Box[\xi](x,y)|^p}{\prod_{k=1}^d |x_k - y_k|^{\alpha_k p + 1}}  dx dy \right]^{1/p} \eqno(4.3)
$$
or equally

$$
|\Box[\xi](x,y)|^p \le  8^{dp} \ 4^{d} \ \prod_{i=1}^d \left[\frac{\alpha_i+1/p}{\alpha_i-1/p} \right]^p \cdot
\vec{\delta}^{\vec{\alpha} p} \cdot  \vec{\delta}^{- 1} \times
$$

$$
 \left[ \int \int \ldots \int_{[0,1]^d} \frac{|\Box[\xi](x,y)|^p}{\prod_{k=1}^d |x_k - y_k|^{\alpha_k p + 1}}  dx dy \right]. \eqno(4.4)
$$
 We get taking expectation:

$$
|\Box[\xi],\vec{\delta}|_p \le \vec{\delta}^{\vec{\alpha}} \cdot \vec{\delta}^{-1/p} \cdot \theta_{\alpha}(p).\eqno(4.5)
$$
 We intend to take the infinum of bide-side inequalities (4.5) over $ p; \ p \in (A,B). $ Note that
$$
\inf_{p \in (A,B)} |\Box[\xi],\vec{\delta}|_p = |\Box[\xi],\vec{\delta}|_A
$$
(Lyapunov's inequality)  and

$$
\inf_{p \in (A,B)} \vec{\delta}^{\vec{\alpha}} \cdot \vec{\delta}^{-1/p} \cdot \theta_{\alpha}(p) =
\vec{\delta}^{\vec{\alpha}} \cdot \frac{1}{\sup_{p \in (A,B)} (\prod_{k=1}^d \delta_k)^{1/p} \theta_{\alpha}(p) }=
\frac{\vec{\delta}^{\vec{\alpha}}}{\phi(G\theta_{\alpha}, \prod_{k=1}^d \delta_k) }.\eqno(4.6)
$$
Q.E.D. \par

\vspace{4mm}

{\bf Remark 4.1.} We can obtain the  exponential bounds for the tail of distribution of r.v. $ \Omega[\xi],\vec{\delta} $
as follows. Let us define the so-called {\it truncated } fundamental function

$$
\phi_q(G\psi, \delta) = \sup_{p \in (q,B)} \frac{\delta^{1/p}}{\psi(p)}, \ A < q < B.
$$
 Denote also

 $$
 \lambda(q, \delta) = \frac{1}{\phi_q(G\theta_{\vec{\alpha}}, \delta)}.
 $$
 It follows from the proposition of theorem 4.1 that

 $$
 \frac{\Omega[\xi],(\vec{\delta})}{\vec{\delta}^{\vec{\alpha}}} \le \lambda(q, \prod_{k=1}^d \delta_k),
 $$
or equally

$$
\left| \left|  \frac{\Omega[\xi],(\vec{\delta})}{\vec{\delta}^{\vec{\alpha}}} \right| \right| G \lambda(\cdot, \prod_{k=1}^d \delta_k)
\le 1.
$$
 It remains to use the conclusion of remark 1.4. \par

\vspace{4mm}

{\bf Theorem 4.2.} {\it  Suppose that for some finite positive constants  }  $  K, \ \alpha, \  \{\beta_k \}, \ k=1,2,\ldots, d $

$$
{\bf E} |\Box [\xi](x,y) |^{\alpha} \le K \cdot \prod_{k=1}^d |x_k - y_k|^{1 + \beta_k}. \eqno(4.7)
$$
{\it Then there exists a non-negative random variable} $  \tau $ {\it with  finite moment of order }
$  \alpha: \ {\bf E} \tau^{\alpha} \le 1 $ {\it such that}

$$
\Omega[\xi](\vec{\delta}) \le C(\alpha, \vec{\beta}, d) \cdot \tau \cdot K^{1/\alpha} \cdot
\prod_{k=1}^d \left[ \delta_k^{\beta_k/\alpha} \ |\log \delta_k|^{1/\alpha} \right], \ 0 < \delta_k  \le 1/e. \eqno(4.8)
$$
{\bf Proof.}  We can assume $  K = 1 $ and use the multidimensional  Garsia-Rodemich-Rumsey inequality, in which
we choose

$$
\Phi(x) = |x|^{\alpha}, \ p_k(x) = |x|^{\gamma_k}, \ 2/\alpha < \gamma_k < (2+\beta_k)/\alpha, \ \gamma_k =  (2+\beta_k)/\alpha - \epsilon_k.
\eqno(4.9)
$$
 Let us introduce the following random variable (r.v.)

$$
B = \int \int \ldots \int_{[0,1]^{2d}} \left| \frac{|\Box[\xi](x,y)|}{\prod_{k=1}^d |x_k - y_k|^{\gamma_k} }  \right|^{\alpha} dx \ dy. \eqno(4.10)
$$
 We have using polar coordinates:

$$
{\bf E} B \le   \int \int \ldots \int_{[0,1]^{2d}} \frac{ \prod_{k=1}^d |x_k - y_k|^{1 + \beta_k}}{ \prod_{k=1}^d |x_k - y_k|^{2 + \beta_k - \gamma_k} } dx \ dy \le
$$

$$
C_2(\alpha, \vec{\beta}, d) \prod_{k=1}^d \int_0^{\sqrt{d}} r_k^{-1+\alpha \epsilon_k}  \ dr_k = C_3(\alpha, \vec{\beta}, d)
 \prod_{k=1}^d (1/\epsilon_k). \eqno(4.11)
$$
 Therefore the r.v. $  B  $ may be represented as a product
 $$
 B =  C_3(\alpha, \vec{\beta}, d) \tau^{1/\alpha} \prod_{k=1}^d (1/\epsilon_k),\eqno(4.12)
 $$
where $ {\bf E} \tau^{\alpha} \le 1. $ \par
 We get substituting into  the multidimensional  Garsia-Rodemich-Rumsey inequality

$$
|\Box \xi|(\vec{\delta}) \le C_4(\alpha, \vec{\beta}, d) \int_0^{\delta_1} \int_0^{\delta_2} \ldots \int_0^{\delta_d}
\left[ \frac{B}{\prod_{k=1}^d u_k^2} \right]^{1/\alpha} \ \left[ \prod_{k=1}^d  \gamma_k u_k^{\gamma_k-1} \ d u_k  \right] =
$$

$$
 C_5(\alpha, \vec{\beta}, d)  \ K^{1/\alpha} \ \tau \ \prod_{k=1}^d \delta_k^{\beta_k/\alpha} \
\prod_{k=1}^d \left[\epsilon_k^{-1/\alpha} \ \delta_k^{-\epsilon_k}  \right]. \eqno(4.13)
$$
 Choosing  $ \epsilon_k = C_6(\alpha, \beta_k, k)/|\log \delta_k|,   $ we arrive to the assertion of theorem  4.2.\par

{\bf Remark 4.2.} Let us show the exactness of assertion of theorem 4.2. It is sufficient to consider simple example.
Let $ d = 1 $ and let $ \xi(t) = w(t), \ t \in [0,1] $ be ordinary Brownian motion. We can choose
$  \alpha = \alpha(\Delta): = 2+2\Delta, \ \beta = \beta(\Delta):= \Delta, \Delta = \const \ge 1.$ Indeed:

$$
{\bf E} |w(t) - w(s)|^{2+2\Delta} = C(\Delta) \ |t-s|^{1+\Delta}, \ s,t \in [0,1],
$$
 Note that

 $$
 \lim_{\Delta \to \infty} \frac{\beta(\Delta)}{\alpha(\Delta)} = 1/2.
 $$

But it is well known that

$$
\overline{\lim}_{\delta \to 0+} \frac{\omega(w,\delta)}{\delta^{1/2} |\log \delta|^{1/2}} > 0
$$
almost everywhere. \par

\vspace{3mm}

  Obtained in this section results specify and generalize  ones  in the articles \cite{Garsia1}, \cite{Hu1}, \cite{Ral'chenko1}. \par
 Another approach to the problem of (ordinary) continuity  of random fields based on the so-called generic chaining method
 and entropy  technique with described applications see in  \cite{Bednorz1}, \cite{Fernique1}, \cite{Kozachenko1}, \cite{Ledoux1},
 \cite{Ostrovsky1}, \cite{Talagrand1}, \cite{Talagrand2} etc.\par

\vspace{3mm}

\section{Concluding remarks.}

\vspace{3mm}

{\bf A. Exactness of our estimations: one dimensional case.} \par

\vspace{3mm}

 Note that at $ \alpha = d  $ and following $ \alpha = d  = 1 $ our estimation in the second section and
the integer order results from the article  \cite{Ostrovsky100} coincides  up to multiplicative constants.
But it is proved in \cite{Ostrovsky100} that the integer order estimates are asymptotically exact.\par

\vspace{3mm}

{\bf B. Exactness of our estimations: multidimensional case.} \par

\vspace{3mm}

 Let us prove that the exponents $\vec{\alpha} - 1/p $ in (3.2) - (3.3) are unimprovable.
Recall that $ p > \max_k (1/\alpha_k). $ As before, $ \psi_{\alpha} (\cdot) = \psi_{(p)}(\cdot). $ \par
 Namely, let us denote in the one-dimensional case $ d = 1 $

$$
V_{\alpha}(f,\delta) =  |\log  \omega(f, \delta)|:
\left|\log \left[\frac{\delta^{\alpha}}{\phi(G\psi_{(p)}, \delta)} \cdot  ||f||G \psi_{(p)} \right] \right|,\eqno(5.1)
$$

$$
\underline{V}_{\alpha}  = \inf_{f \in G\psi_{(p)}(\cdot), f \ne \const} \overline{\lim}_{\delta \to 0+} V_{\alpha}(f,\delta). \eqno(5.2)
$$
It follows from theorem 2.1 that  $ \underline{V}_{\alpha}  \ge 1; $ let us prove
the opposite inequality. Consider the following example (more precisely, the family of examples):

$$
f_{\Delta, \alpha}(x) = x^{\alpha - 1/p + \Delta}, \ \Delta = \const \in (0, 1-\alpha  + 1/p). \eqno(5.3)
$$
 Obviously,
 $$
 \omega(f_{\Delta, \alpha},\delta)=  \delta^{\alpha - 1/p + \Delta}.
 $$
  Further, it is no hard to compute using polar coordinates: as $ \delta \to 0+ $

 $$
\frac{\delta^{\alpha}}{\phi(G\psi_{(p)}, \delta)} \cdot  ||f_{\Delta, \alpha}||G \psi_{(p)}  \sim  C(\alpha,p, \Delta) \
\delta^{\alpha - 1/p}.
 $$
 Since the value $ \Delta $ is arbitrary, we conclude  $ \underline{V}_{\alpha} \le 1.$\par

 The multidimensional example may be constructed as a factorable function of a view

 $$
 g_{\Delta, \vec{\alpha}}(\vec{x}) = \prod_{k=1}^d f_{\Delta, \alpha_k}(x_k). \eqno(5.4)
 $$

\vspace{3mm}

{\bf C. Simplification of our multidimensional estimate.} \par

\vspace{3mm}
 Let us consider in this  paragraph the following important coefficient:

 $$
 L = \prod_{k=1}^d \frac{\alpha_k + 1/p}{\alpha_k-1/p}, \ p > 1/\alpha_0, \ \alpha_0 := \min_k \alpha_k,
 $$
meaning to extract the main factor as $ p \to 1/\alpha_0. $ \par
 We have:

$$
L = \prod_{k: \alpha_k > \alpha_0} \frac{\alpha_k + 1/p}{\alpha_k-1/p} \times
\prod_{k: \alpha_k = \alpha_0} \frac{\alpha_k + 1/p}{\alpha_k-1/p} \le
$$

$$
\prod_{k: \alpha_k > \alpha_0} \frac{\alpha_k + \alpha_0}{\alpha_k-\alpha_0} \times
\left\{ \frac{\alpha_0 + 1/p}{ \alpha_0-1/p} \right\}^M.
$$

\vspace{3mm}

{\bf D. General rectangle distance.} \par

\vspace{3mm}

{\ 1.} Let$ X_j = \{x_j\}, j = 1,2,\ldots,d $ be arbitrary sets and $ f:  Z= \otimes_{j=1}^d X_j \to R $ be numerical function.
Define the following function

$$
\rho_f(\vec{x}, \vec{y}) =  \rho_f(x, y) = |\Box[f](x,y)|, \ x,y \in Z.
$$

 Note the following properties of the function $ \rho_f(x, y). $
 $$
 \leqno(a) \hspace{3mm}  \rho_f(x,y) \ge 0; \ \exists j = 1,2,\ldots,d  \hspace{3mm}  x_j = y_j
 \Rightarrow \rho_f(x, y) = 0.
 $$
 (non-negativity);
 $$
\leqno(b) \hspace{3mm}   \rho_f(x, y) = \rho_f(y, x),
 $$
 (symmetry);
 $$
\leqno(c)  \hspace{3mm}  \rho_f(x, z) \le \rho_f(x,y) + \rho_f(y,z), \ x,y,z \in Z,
 $$
 (rectangle inequality).\par

\vspace{3mm}

{\bf Definition of  a rectangle  distance.} \par
 Arbitrary numerical function of $ 2 d $ variables $ \rho(x,y), \ x,y \in Z  $ which satisfies the properties
 (a,b,c) is said to be a {\it rectangle  distance.} \par

 \vspace{3mm}

{\bf Example.} Let $  \xi = \xi(x), \ x \in Z $ be a random field with condition
$$
\exists q \in [1,\infty], \  \sup_{x \in Z} \left[ {\bf E} |\xi(x)|^q  \right]^{1/q} < \infty.
$$
  The function

 $$
 \rho^{(\xi,q)}(x,y) =   \left\{ {\bf E} |\Box[\xi](x,y)|^q  \right\}^{1/q}
 $$
 is bounded {\it natural} rectangle  distance generated by r.f. $ \xi = \xi(x),  $  regarded before. \par
  Obviously, instead the classical $ L_q $ norm may be used arbitrary rearrangement norm, for
 instance, Orlicz, Grand Lebesgue, Lorentz or Marcinkiewicz norm etc. \par

\vspace{3mm}

{\bf E. Scaling method.} \par

\vspace{3mm}

  We intend to prove here the exactness of inequality (1.5) by means of the so-called scaling method, see
 \cite{Stein1}, \cite{Talenti1}.  Indeed, we can extrapolate the function $ f $ in (1.5)  as continuous function
in the closed interval $ [0,2] $ with support  $ [0,2] $ such that on the set $  [1,2] \ f $   is linear, $ f(2) = 0.  $
For such a function (1.5) also holds.\par
Introduce the dilation operator $ T_{\lambda}[f] = f(\lambda x), \ \lambda > 0. $  Consider the following strengthening of  (1.5)
for any continuous function with compact support belonging to the space $ C^{(0)}_{\alpha - 1/p}, $ where by definition the
space $ C^{(0)}_{\beta},  \beta \in (0,1]$ consists on all (continuous) function with finite semi-norm

$$
||f||C^{(0)}_{\beta} = \sup_{\delta > 0}  \frac{\omega(f,\delta)}{\delta^{\beta}}:
$$

$$
\omega(f,\delta) \le 8 \cdot 4^{1/p} \cdot \left[\frac{\alpha + 1/p}{\alpha - 1/p}\right]   \cdot \delta^{\alpha - 1/p} \cdot \gamma(\delta) \cdot
\left[ \int_0^{\infty} \int_0^{\infty} \frac{|f(x) - f(y)|^p \ dx dy }{|x-y|^{\alpha p + 1}}   \right]^{1/p} =:
$$
$$
 \delta^{\alpha - 1/p} \cdot \gamma(\delta) ||f||U(\alpha,p),    \eqno(5.5)
$$
where $ \lim_{\delta \to 0+} \gamma(\delta) = 0. $\par
 Applying (5.5) for the non-constant function $  T_{\lambda}[f], $ we obtain after simple calculations:

$$
||T_{\lambda}f||^p U(\alpha,p) = \int_0^{\infty} \int_0^{\infty}  \frac{|f(\lambda x) - f(\lambda y)  |^p \ dx \ dy }{|x-y|^{\alpha p + 1}}=
$$

$$
\lambda^{-2} \ \int_0^{\infty} \int_0^{\infty} \frac{|f(x) - f(y)|^p \ dx \ dy}{ |x/\lambda - y/\lambda |^{\alpha p + 1}} =
\lambda^{-1 + \alpha p }\int_0^{\infty} \int_0^{\infty}  \frac{|f(x) - f(y) |^p \ dx \ dy }{|x-y|^{\alpha p + 1}}=
$$

$$
\lambda^{-1 + \alpha p } ||f||^p U(\alpha,p), \ ||T_{\lambda}f|| U(\alpha,p) = \lambda^{\alpha - 1/p} ||f|| U(\alpha,p);
$$

 $$
\frac{\omega(f, \lambda \delta)}{ (\lambda \delta)^{\alpha - 1/p} } \le  \gamma(\delta) ||f||U(\alpha,p). \eqno(5.6)
 $$
We get taking  supremum over $ \lambda: $

$$
||f||C^{(0)}_{\alpha - 1/p}  \le  \gamma(\delta) ||f||U(\alpha,p),
$$
 which is not true as $ \delta \to 0+. $ \par

\vspace{3mm}

{\bf F. General Orlicz approach.} \par

\vspace{3mm}

 Let $ \Phi = \Phi(u) $ be again the Young-Orlicz function. We will denote the Orlicz norm by means of the function $  \Phi $
 of a r.v. $ \kappa $ defined on our probabilistic space  as $ |||\kappa|||L(\Phi).   $ \par
  We introduce  the natural rectangle distance $ \rho_{\Phi}(x,y) $ as follows:

  $$
  \rho_{\Phi}(x,y):= |||\Box[\xi](x,y), \ x,y \in D = [0,1]^d, \eqno(5.7)
  $$
so that  for the r.v.

$$
Y = \int_D \int_D \Phi \left( \frac{\Box[\xi](x,y)}{\rho_{\Phi}(x,y}  \right) \ dx \ dy
$$
we have

$$
 {\bf E} Y  = \int_D \int_D {\bf E} \Phi \left( \frac{\Box[\xi](x,y)}{\rho_{\Phi}(x,y)}  \right) \ dx \ dy \le 1, \eqno(5.8)
$$
since $ \int_D \int_D  dx \ dy = 1.  $\par
 Let also $ \rho^{(\Phi)}(x-y) $ be translation invariant strictly increasing continuous distance majored  $ \rho_{\Phi}(x,y):$

$$
\rho_{\Phi}(x,y) \le \rho^{(\Phi)}(x-y).
$$
 We denote the particular distances

$$
p_k(|y_k - x_k|) =  \rho^{(\Phi)} (1,1,  \ldots,1, |x_k-y_k|, 1, \ldots,1).
$$
 It follows immediately from the multidimensional version of Garsia-Rodemich-Rumsey inequality that

 $$
 \Omega[\xi](\vec{\delta}) \le  8^d \int_0^{\delta_1} \int_0^{\delta_2} \ldots \int_0^{\delta_d}
 \Phi^{-1} \left( \frac{4^d Y}{\prod_{k=1}^d u_k^2} \right) \prod_{k=1}^d d p_k(u_k).\eqno(5.9)
 $$

  Of course, the inequality (5.9) is pithy if the integral the right-hand side  convergent;
  then the r.f. $  \xi(\cdot) $ is continuous with probability one. \par
   The Gaussian (more precisely, subgaussian) case considered in \cite{Garsia1}, \cite{Hu1}, \cite{Ral'chenko1}
 may be obtained by choosing $ \Phi(z) = \exp(z^2/2) - 1. $ It may be considered easily the example when
 $ \Phi(z) = \exp(|z|^m/m) - 1, \ m = \const > 0. $  \par

\vspace{3mm}

{\bf G.  Fractional Orlicz - Sobolev inequalities. }

\vspace{3mm}

 Let $ f:D = [0,1]^d \to R  $ be (measurable) function. We define the following natural $  \Psi $ function depending
 on the vector positive parameter $ \vec{\alpha}:  $

$$
\tau_{\vec{\alpha}}(p) = || G_{\vec{\alpha}}[f]{\cdot,\cdot} ||_{p,D^2,\nu}, \ p > 1/\min(\alpha_k). \eqno(5.10)
$$
 It will be presumed that the function $ \tau_{\vec{\alpha}}(p) $ there exists:

 $$
 \supp \tau_{\vec{\alpha}}(\cdot) = (A,\infty), A > 1/\min \alpha_k.
 $$
 We can construct the following exponential $  N = N_{\vec{\alpha}}$ Young-Orlicz function as in remark 1.3:

 $$
 N_{\vec{\alpha}}(u) = e^{[p \log \tau_{\vec{\alpha}}(p)]^*(\log |u|)}, \ |u|  > 3. \eqno(5.11)
 $$
 {\it We offer in this subsection  a multidimensional (rectangle) version  of fractional Orlicz-Sobolev inequality for the exponential Orlicz's space}
 $ L ( N_{\vec{\alpha}}(\cdot)).  $ Note that the integer ordinary (interval)  Orlicz-Sobolev inequality for the arbitrary Orlicz's space
 is considered, e.g. in \cite{Adams1}, chapter 11; \cite{Rao1}, chapter 9. \par
 We  infer on the basis of theorem 3.1 and remark 1.3:

 {\bf Proposition 5.1.}

 $$
 \Omega[f](\vec{\delta}) \le  C(\vec{\alpha}, d) \ \frac{\vec{\delta}^{\vec{\alpha}} \cdot ||f||L( N_{\vec{\alpha}}(\cdot))}{\phi(G\tau_{\alpha}, \prod_{k=1}^d \delta_k )}. \eqno(5.12)
 $$
 As a slight strengthening: \\
{ \bf Proposition 5.2.}

 $$
 \Omega[f](\vec{\delta}) \le \inf_{\vec{\alpha}}
 \left[C(\vec{\alpha}, d) \ \frac{\vec{\delta}^{\vec{\alpha}} \cdot ||f|| L( N_{\vec{\alpha}}(\cdot))}{\phi(G\tau_{\alpha}, \prod_{k=1}^d \delta_k )} \right]. \eqno(5.13)
 $$

   Note in addition that the fundamental function $ \phi(L(\Phi), \delta)  $ for arbitrary probabilistic  Orlicz's
spaces $ L(\Phi)  $ is calculated, e.g.  in  the classical book of Krasnoselsky M.A., Rutizky Ya.B.  \cite{Krasnoselsky1},
chapter 2, section 9:

$$
\phi(L(\Phi), \delta)  = \delta \cdot \Phi^{-1}(1/\delta). \eqno(5.14)
$$

 See also more modern books \cite{Rao1}, \cite{Rao2}.\par

\vspace{9mm}

{\bf   Acknowledgements.} The authors  would very like to thank  prof. S.V. Astashkin, M.M.Milman, and  L.Maligranda for sending
of Yours remarkable papers and several useful suggestions. \par

\end{document}